\documentclass[12pt]{elsarticle}

\usepackage{amssymb}
\usepackage{amsmath}
\usepackage{mathrsfs}
\usepackage{amsfonts}

\newtheorem{example}{Example}[section]
\newtheorem{definition}{Definition}[section]
\newtheorem{lemma}{Lemma}[section]
\newtheorem{proposition}{Proposition}[section]
\newtheorem{theorem}{Theorem}[section]
\newtheorem{remark}{Remark}[section]

\usepackage{comment}
\usepackage{color}

\newenvironment{proof}[1][\noindent \textit{Proof. }]{#1}{ \hfill $\square$ \vspace{2mm}}

\journal{?}

\begin{document}

\begin{frontmatter}

\title{Globally solvable complexes of pseudo-differential operators on the torus}

\author{\corref{cor1}Fernando de Ávila Silva} \ead{fernando.avila@ufpr.br}

\author{Cleber de Medeira} 
\ead{clebermedeira@ufpr.br}

% Author affiliation
\affiliation{organization={Departamento de Matemática  UFPR},%Department and Organization
            addressline={Caixa Postal 19081}, 
            city={Curitiba},
            postcode={CEP 81531-980}, 
            state={Paraná},
            country={Brazil}}

%% Abstract
\begin{abstract}
We consider a complex of pseudo-differential operators  associated  with an overdetermined system of operators defined on the torus. We characterize  the global solvability of this complex when the system has constant coefficients. Furthermore, we
establish some connections between the global solvability of the complex  and the global solvability of its normal form.	
\end{abstract}

%% Keywords
\begin{keyword}
Global solvability\sep complex of pseudo-differential  operators \sep 
overdetermined systems \sep Liouville forms, Fourier series.
\MSC[2020] 35A01  \sep 35N10 \sep 58J10
\end{keyword}

\end{frontmatter}

\section{Introduction}

This work deals with the global solvability for a class of systems of pseudo-differential operators on the torus $\mathbb{T}^{n+N}\simeq \mathbb{R}^{n+N}/2\pi
\mathbb{Z}^{n+N}$, at the level of $p$-forms, with $p=0,1,\ldots,n-1$. The periodic coordinates in $\mathbb{T}^{n+N}$ are written as $(t, x)=(t_1, \dots, t_n,x_1,\ldots,x_N)$ and the system is given by
\begin{equation}\label{main system}
	L_j=D_{t_j}+c_j(t)P_j(D_x),\quad j=1,\ldots,n,
\end{equation}
where $D_{t_j}=i^{-1}\partial_{t_j}$,  $c_j$ is a smooth function on $\mathbb{T}_t^n$ and $P_j(D_x)$ is a pseudo-differential operator of order $m_j\in\mathbb{R}$ on $\mathbb{T}_x^N$,
defined by 
\begin{equation*}\label{pdo}
	P_j(D_x)  u(x) = \sum_{\xi \in \mathbb{Z}^N}{e^{i x \cdot \xi} p_j(\xi) \widehat{u}(\xi)},
\end{equation*}
where  
$
\widehat{u}(\xi) = (2\pi)^{-N} \int_{\mathbb{T}^N}{e^{- i x \cdot \xi} u(x) dx}
$  denotes the Fourier coefficients of $u$ and  the toroidal symbol  $p_j(\xi)$ belongs to the class $S^{m_j}(\mathbb{Z}^N)$ in the sense of Definition 4.1.7 in \cite{RT3}. In particular, there exists $C>0$ such that
\begin{equation}\label{symbol pj}
	|p_j(\xi)|\leq C\langle\xi\rangle^{m_j},
\end{equation}
for all $j=1,\ldots,n$ and $\xi\in\mathbb{Z}^N$. Also, 
\begin{equation*}\label{recover_symbol}
	p_j(\xi) = e^{-i x\cdot \xi} P_j(D_x) e^{i x \cdot \xi}.
\end{equation*}

\bigskip

We assume that  $c(t,\xi)\doteq\sum_{j=1}^{n}p_j(\xi)c_j(t)dt_j$ is a  closed $1$-form, for each $\xi\in\mathbb{Z}^N$. Therefore,  in a natural way, we may  associate System \eqref{main system} with a complex of  pseudo-differential operators, according with Treve's seminal work \cite{Tre76}. We now move on to describing it.

For each $p=0,1,\ldots,n-1$, we  denote by $C^\infty(\mathbb
T^{n+N};\wedge^{p,0})$ the space of  $p$-forms   on  $\mathbb{T}^{n+N}$ spanned by   $dt_K=dt_{k_1}\wedge\ldots \wedge dt_{k_p}$, where $K=(k_1,\ldots,k_p)$ is a multi-index  of positive integers such that $1\leqslant k_1<\cdots<k_p\leqslant n$ and $|K|$ is the length of $K$. Thereby, 
each $u\in  C^\infty(\mathbb
T^{n+N};\wedge^{p,0})$ can be written as 
$$
u(t,x)=\sum_{|K|=p}u_K(t,x)dt_K,
$$
where $u_K\in C^{\infty}(\mathbb{T}^{n+N})$. In an equivalent way, we denote by $\mathcal{D}'(\mathbb{T}^{n+N};\wedge^{p,0})$ the space of $p$-currents on $\mathbb{T}^{n+N}$ spanned by $dt_K=dt_{k_1}\wedge\ldots \wedge dt_{k_p}$, with coefficients in $\mathcal{D}'(\mathbb{T}^{n+N})$.

Thus, associated  with System \eqref{main system}, we consider the operator
$$
\mathbb{L}^p:\mathcal{D}'(\mathbb{T}^{n+N}; \wedge^{p,0})\rightarrow \mathcal{D}'(\mathbb{T}^{n+N};\wedge^{p+1,0})
$$
defined by 
$$
\mathbb{L}^p u=\sum_{j=1}^{n}\sum_{|K|=p}(iL_ju_K)dt_j\wedge dt_K,
$$
where   $u=\sum_{|K|=p}u_K(t,x)dt_K$.

\bigskip

Now, consider the operator $c(t,D_x)$ acting in $u$ as follows
$$
c(t,D_x)\wedge u\doteq \sum_{j=1}^{n}\sum_{K|=p|}ic_j(t)P_j(D_x)u_K\;dt_j\wedge dt_K.
$$
Therefore, we may rewrite the operator $\mathbb{L}^p$ as
\begin{equation}\label{main operator}
	\mathbb{L}^p=d_t+c(t,D_x)\wedge,
\end{equation}
where $d_t$ stands for the exterior derivative on $\mathbb{T}_t^n$.

The assumption that $c(t,\xi)$ is a closed 1-form, for every $\xi\in \mathbb{Z}^N$, implies that $[L_j,L_k]=0$, for all $j,k$.   Then,  if  $f\in C^\infty(\mathbb
T^{n+N};\wedge^{{p+1},0})$ satisfies $\mathbb{L}^{p}u=f$, for some $u=\sum_{|K|=p}u_K(t,x)dt_K  \in\mathcal{D}'(\mathbb{T}^{n+N}, \wedge^{p,0})$,
we get
\begin{equation}\label{first comp condition}
	\mathbb{L}^{p+1}f=\mathbb{L}^{p+1}(\mathbb{L}^pu)=\sum_{|K|=p}\sum_{1\leq j<k\leq n }i[L_j,L_k](u_K)dt_j\wedge dt_k\wedge dt_K=0,
\end{equation}
concluding that $Im(\mathbb{L}^p)\subset \ker(\mathbb{L}^{p+1})$.
This means that $\mathbb L^p$  defines a complex of pseudo-differential operators on $\mathbb{T}^{n+N}$ as described in Section I.3. in  \cite{Tre76}.

Hence, in this article, our aim is to carry out a study of the global sol\-va\-bi\-lity of $\mathbb{L}^p$ at each level of this complex of operators.

\bigskip

We stress that the global solvability of linear differential operators of type \eqref{main operator}, associated with closed 1-forms and defined on smooth compact manifolds, has been studied in different scales of functional spaces.  We highlight the case $\mathbb{L}^p=d_t+c(t)\wedge\partial_x$  which is widely explored by many authors, including, for instance, 
Cardoso and Hounie \cite{CarHou77},
Bergamasco and Petronilho \cite{BerPet99jmaa}, 
Ber\-ga\-masco et al. \cite{BerKirNunZan12}, Bergamasco, Medeira and Zani \cite{BerMedZan12}, Zugliani and Hounie \cite{HouZug2017,HouZug2019}. More recently, still considering the same scenario, results in Gevrey framework also became of interest as we may see in Araujo, Dattori and Lessa \cite{ADL23}, Bergamasco, Medeira and Zani \cite{BerMedZan21} and Dattori and Meziani \cite{Silva2020}. 
Regarding the global hypoellipticity of the operator $\mathbb{L}^p$, at the first level of the complex ($p=0$), in the $C^{\infty}$ case, we call attention to the authors Bergamasco, Cordaro and Malaguti  \cite{BerCorMal93}; for the  analytic hypoellipticity we mention  Bergamasco \cite{Ber99}; lastly, for Gevrey classes, we cite Arias, Kirilov and Medeira \cite{AriKirMed19}.

However, less common are investigations for  classes of pseudo-differential operators. For the global hypoellipticity problem, at the first level of the complex, we can cite  
Ávila and Medeira \cite{AviMed21}, and Ávila et al. \cite{AviGonKirMed19}.  Also, the machinery of pseudo-differential operators can be applied to study certain perturbations of differential operators, as considered by Ferra, Petronilho and Lessa \cite{FerraPetrLessa}, Braun et al.  \cite{BraunChinniCordaroJah} and Chinni and Cordaro  \cite{ChinniCordaro2017}.  Therefore, as far as we know, our work presents the first results for the global solvability problem of the complex of pseudo-differential operators $\mathbb{L}^p$ on the torus.

\bigskip

The paper is structured as follows: In Section \ref{comp cond} we begin the discussion of our question establishing some  compatibility conditions for the existence of a solution $u$ to the equation $\mathbb{L}^pu =f\in C^\infty(\mathbb T^{n+N};\wedge^{p+1, 0})$. Thus, we can precisely define the sense of global solvability that we are interested in studying.

In Section \ref{section coefficients constants }, by assuming constant coefficients in System \eqref{main system}, we  completely  characterize the global solvability of the associated operator  $\mathbb{L}^p$ in terms of the asymptotic behavior of its symbol. The main result is Theorem \ref{constant coefficients thm} and as an illustration of this result, in Subsection \ref{homo-symbols}, we present a characterization of the global solvability for a class of operators with homogeneous symbols. In this case, the result obtained in Theorem \ref{homogeneous operators} is closely related with certain Diophantine approximations of irrational numbers by rational numbers with common denominators, what includes, for instance, the well-known Liouville numbers (vectors).

Section \ref{section reduction normal form} is dedicated to state some connections between the global sol\-va\-bi\-lity of the operator $\mathbb{L}
^p$ and the global solvability of the associated operator with constant coefficients \eqref{normal form operator}, which is named normal form. We state a condition (see \eqref{condition D}) that guarantees the reduction of the study of global solvability to the normal form, as stated in Theorem \ref{reduction normal form thm}. Next, we present some classes of operators that satisfy this reduction condition. Also, we include some examples in order to better understand the results obtained. In particular, in Subsection \ref{sec_dec_system}, we illustrate our results considering a class of decoupled systems.

\section{Compatibility conditions}\label{comp cond}

There are natural compatibility conditions   on 
$f\in C^\infty(\mathbb{T}^{n+N};\wedge^{{p+1},0})$ for the existence of a solution $u\in\mathcal{D}'(\mathbb{T}^{n+N}; \wedge^{p,0})$  to the equation $\mathbb{L}^pu=f$. We now move on to describing them. First, we recall by \eqref{first comp condition} that, if there exists a solution $u$ to the equation $\mathbb{L}^pu=f$, then we obtain the first compatibility condition $\mathbb{L}^{p+1}f=0$. 

Now, consider the partial Fourier series of $u=\sum_{|K|=p}u_K(t,x)dt_K$ and $f=\sum_{|J|=p+1}f_J(t,x)dt_J$, with respect to $x$, given by
$$
u(t,x)=\sum_{\xi\in \mathbb{Z}^N}
\widehat{u}(t,\xi)e^{i\xi \cdot x}\quad \textrm{ and }\quad  f(t,x)=\sum_{\xi\in \mathbb{Z}^N}
\widehat{f}(t,\xi)e^{i\xi\cdot  x},
$$
where $\widehat{u}(t,\xi)=\sum_{|K|=p} \widehat{u_K}(t,\xi)dt_K$ and
$\widehat{u_K}(t,\xi)$ is the Fourier transform of $u_K$ with respect to $x$. The coefficients $\widehat{f}(t,\xi)$ are
analogously defined.

For each $\xi\in\mathbb{Z}^N$, we may decompose the closed $1$-form $c(t,\xi)$ as
\begin{equation*}\label{decomposition closed forms}
	c(t,\xi)=c_{\xi0}+d_t \mathscr{C}_\xi(t),
\end{equation*}
where $c_{\xi0}\in\bigwedge ^1\mathbb{C}^n
$
is a constant $1$-form and
$\mathscr{C}_\xi(t)$ is a complex valued smooth function on $\mathbb{T}_t^n$. 
A simple computation shows that
$$
c_{\xi0} = \sum_{j=1}^{n}c_{j0}p_j(\xi) dt_j,
$$
where  $$
c_{j0} = \dfrac{1}{2\pi}\int_{\mathbb{T}^1}c_j(0, \ldots, t_j, \ldots, 0)dt_j.
$$

In particular, when  $c_{\xi0}$ is a real  1-form, as in Definition 2.1 in \cite{BerCorMal93},  we say that $c_{\xi0}$  is integral if
$$
\left(\frac{1}{2\pi}\int_{\sigma}c_{\xi0} \right)\in\mathbb{Z},
$$
for any 1-cycle $\sigma$ in $\mathbb{T}^n$.  Furthermore, $c_{\xi0}$ is said to be rational when  there exists a positive integer  $q$ such that $q c_{\xi0}$ is an integral 1-form.

\bigskip

Using the previous definitions, if $\mathbb{L}^pu=f$ we may conclude that
\begin{eqnarray}\label{cond compat}
	\widehat{f}(t,\xi)e^{i(\psi_\xi(t)+\mathscr{C}_\xi(t))} \;\textrm{is a  $(p+1)$-exact form,  whenever $c_{\xi0}$ is integral,\; }
\end{eqnarray}
where   $\psi_{\xi}\in{C}^{\infty}(\mathbb{R}^n;\mathbb{R})$ is such that $d\psi_\xi=\Pi^*( c_{\xi0})$  and $\Pi:\mathbb{R}^n\rightarrow \mathbb{T}^n$ denotes the universal covering space of $\mathbb{T}^n$. 

In fact, replacing the Fourier series in the equation $\mathbb{L}^pu=f$, we obtain
$$
d_t\widehat{u}(t,\xi)+ic(t,\xi) \wedge\widehat{u}(t,\xi)=\widehat{f}(t,\xi),\quad \xi\in\mathbb{Z}^N,
$$
or equivalently, 
$$
(d_t+i c_{\xi0}\wedge)e^{i\mathscr{C}_\xi(t)}\widehat{u}(t,\xi)=e^{i\mathscr{C}_\xi(t)}\widehat{f}(t,\xi),\quad \xi\in\mathbb{Z}^N.
$$

In particular, if $\xi\in\mathbb{Z}^N$ is such that $c_{\xi0}$ is integral  and $\psi_{\xi}\in{C}^{\infty}(\mathbb{R}^n;\mathbb{R})$ is such that $d\psi_\xi=\Pi^*( c_{\xi0})$, then by Lemma 2.3 in \cite{BerCorMal93}, we have
$$
P,Q\in \mathbb{R}^n, \quad \Pi(P)=\Pi(Q)\Rightarrow (\psi_\xi(P)-\psi_\xi(Q))\in 2\pi \mathbb{Z}.
$$	
Thus, $e^{i\psi_\xi(t)}$ is a well defined function  on $\mathbb{T}_t^n$. Also,
$$
d_t\big(e^{i(\psi_\xi(t)+\mathscr{C}_\xi(t))}\widehat{u}(t,\xi)\big)=\widehat{f}(t,\xi)e^{i(\psi_\xi(t)+\mathscr{C}_\xi(t))},
$$
implying that $e^{i(\psi_\xi(t)+\mathscr{C}_\xi(t))}\widehat{f}(t,\xi)$ is exact.

Motivated by the previous arguments, for each $p=0,1,\ldots,n-1$, we define the following compatibility conditions set:
\begin{center}
	$\mathbb{E}^p=\big\{f\in C^\infty(\mathbb T^{n+N};\wedge^{ {p+1},0});\;\mathbb{L}^{p+1}  f=0\textrm{ and $(\ref{cond compat})$ holds} \big\}$.
\end{center}

\begin{definition}
	We say that the operator $\mathbb{L}^p$ is globally solvable  on $\mathbb{T}^{n+N}$ if for each $f\in\mathbb{E}^p$
	there exists $u\in \mathcal{D}'(\mathbb{T}^{n+N};\wedge^{p,0})$ such that $\mathbb{L}^p u=f.$
\end{definition}

\section{Systems with constant coefficients}\label{section coefficients constants }

In this section we characterize the global solvability of the operator $\mathbb{L}^p$ associated with the following system with constant coefficients defined in $\mathbb{T}^{n+N}$
\begin{equation}\label{system constant coefficients}
	L_j=D_{t_j}+P_j(D_x),\quad j=1,\ldots,n.
\end{equation}

Note that we are incorporating the constant coefficients into the symbols of $P_j(D_x)$. In this case, the associated operator is given by 
\begin{equation}\label{operatorconstant coefficients}
	\mathbb{L}^p=d_t+c(D_x)\wedge,
\end{equation} 
where 
$$
c(D_x)\wedge u=\sum_{j=1}^{n}\sum_{|K|=p}iP_j(D_x)u_K\; dt_j\wedge dt_K.$$
Also, Condition \eqref{cond compat} can be rewritten as follows:  
\begin{eqnarray}\label{comp condition const coef}
	\widehat{f}(t,\xi)e^{i\psi_\xi(t)}\text{ is exact when $c_{\xi0}=\sum_{j=1}^{n}p_j(\xi)dt_j$ is integral}, 
\end{eqnarray} 
where   $\psi_{\xi}\in{C}^{\infty}(\mathbb{R}^n;\mathbb{R})$ is such that $d\psi_\xi=\Pi^*( c_{\xi0})$.\\ 

For each $(\eta,\xi)\in\mathbb{Z}^{n+N}$, we denote by  $\widehat{\mathbb{L}}(\eta,\xi)\doteq i\sum_{j=1}^n(\eta_j+p_j(\xi))dt_j$ the symbol of $\mathbb{L}^p$ and we consider its norm given by
$$
\|\widehat{\mathbb{L}}(\eta,\xi)\| = \max_{j=1, \ldots, n} |\eta_j + p_j(\xi)|. 
$$

Before we state  the main results of this section, we recall the following key lemma, whose proof can be found  in \cite{BerPet99jmaa}, Lemma 2.1. 

\begin{lemma}\label{lemma solution constant coef system}
	Let $\mathscr{L}=\sum_{j=1}^{n}\mathscr{L}_jdt_j \in \bigwedge^{1}\mathbb{C}^n$ be a nonzero constant 1-form. Given $\mathscr{F}=\sum_{|J|=p+1}\mathscr{F}_J dt_J$  in $ \bigwedge^{p+1}\mathbb{C}^n$, $p=0,1,\ldots,n-1$, the equation 
	\begin{equation}\label{eq Lemma 1}
		\mathscr{L}\wedge \mathscr{U}=\mathscr{F}
	\end{equation}
	has a solution $\mathscr{U}=\sum_{|K|=p}\mathscr{U}_K dt_K$  in   $ \bigwedge^p\mathbb{C}^n$ if and only if 
	\begin{equation}\label{eq Lemma 2}
		\mathscr{L}\wedge \mathscr{F}=0.
	\end{equation} 
	In particular, when \eqref{eq Lemma 2} holds,  a solution of \eqref{eq Lemma 1} is given by 
	\begin{equation}\label{solution lemma}
		\mathscr{U}_0=  \underset{\mu \in J}{\sum_{|J|=p+1}}  (-1)^{\text{sign}(\mu ,j_1,\ldots,\widehat{\mu },\ldots,j_{p+1})} (\mathscr{L}_\mu)^{-1}  \mathscr{F}_J  dt_{J-\{\mu \}},
	\end{equation}
	where 
	$\mathscr{L}_\mu\neq 0 $ and $J-\{\mu\}=(j_1,\ldots,\widehat{\mu },\ldots,j_{p+1})$, whenever $\mu\in J$.
	
	The general solution of \eqref{eq Lemma 1}  is given by 
	$$
	\mathscr{U}= \mathscr{U}_0+\mathscr{L}\wedge \mathscr{W},
	$$
	where $\mathscr{W}$ is an arbitrary element in $\bigwedge^{p-1}\mathbb{C}^n$. Here, we are considering $\bigwedge^{-1}\mathbb{C}^n=\{0\}$. 
\end{lemma}

We consider the set $\mathcal{Z}=\{\xi\in\mathbb{Z}^N;\; c_{\xi0} \text{ is integral}\}$. 
When $\mathcal{Z}\neq \emptyset$, inspired by Section 4 in \cite{BerCorPet96}, we define the space
\begin{equation*}
	\mathcal{D}_{\mathcal{Z}}'(\mathbb{T}^{n+N}) = \left\{ u\in \mathcal{D}'(\mathbb{T}^{n+N});\; u(t,x) = \sum_{\xi\in \mathcal{Z}} \widehat{u}(t,\xi)  e^{i\xi \cdot x} \right\},
\end{equation*}
and we denote by $\mathcal{D}_{\mathcal{Z}}'(\mathbb{T}^{n+N};\wedge^{p,0})$    the space of $p$-currents $u=\sum_{|J|=p}u_Jdt_J$, where  $u_J\in	\mathcal{D}_{\mathcal{Z}}'(\mathbb{T}^{n+N}) $. The space $C_{\mathcal{Z}}^{\infty}(\mathbb{T}^{n+N};\wedge^{p,0})$  can be described in an equivalent way. 

Additionally, we consider   $\mathbb{L}_{\mathcal{Z}}^p$ as  the operator $\mathbb{L}^p$  acting on $\mathcal{D}_{\mathcal{Z}}'(\mathbb{T}^{n+N};\wedge^{p,0})$, that is, $\mathbb{L}_{\mathcal{Z}}^p=\mathbb{L}^p|_{\mathcal{D}_{\mathcal{Z}}'(\mathbb{T}^{n+N};\wedge^{p,0})}$.  In this case, the corresponding compatibility set $\mathbb{E}_{\mathcal{Z}}^p$ and the definition of the global solvability of $\mathbb{L}_{\mathcal{Z}}^p$ are  similar to those made in Section \ref{comp cond}.

\begin{proposition}\label{proposition integral form}
	For each $p=0,1,\ldots,n-1$, the operator $\mathbb{L}_{\mathcal{Z}}^p$ is globally solvable.  
\end{proposition}
\begin{proof}
	Let $f\in \mathbb{E}_{\mathcal{Z}}^p$. If $\xi\in\mathcal{Z}$, then  $c_{\xi0}$ 
	is integral and  by \eqref{comp condition const coef}    the $(p+1)$-form $\widehat{f}(t,\xi)e^{i\psi_\xi(t)}$ is exact. We define $h\in C^{\infty}(\mathbb{T}^{{n+N}};\wedge^{p+1,0})$ by 
	$$
	h(t,x)=\sum_{\xi\in\mathbb{Z}^N}\widehat{h}(t,\xi)e^{i\xi \cdot x},
	$$
	where 
	$$
	\widehat{h}(t,\xi)=
	\begin{cases}
		\widehat{f}(t,\xi)e^{i\psi_\xi(t)},	&\xi\in
		\mathcal{Z}\\
		0, & \xi\in\mathbb{Z}^N\setminus \mathcal{Z}.
	\end{cases}
	$$
	Therefore, $\widehat{h}(t,\xi)$ is exact,
	for all $\xi\in\mathbb{Z}^N$. 
	
	We claim that there exists  a solution $v\in C^{\infty}(\mathbb{T}^{{n+N}};\wedge^{p,0})$ to the equation $d_t v=h$. Indeed, we will present a solution $v$ by defining its Fourier coefficients $\widehat{v}(\eta,\xi)$.

	Since $\widehat{h}(t,\xi)$ is exact, if  $\eta=0$ then $\widehat{h}(\eta,\xi)=0$, for all $\xi\in\mathbb{Z}^N$. In this case, we take  $\widehat{v}(0,\xi)=0$, for all $\xi\in\mathbb{Z}^N$.
	
	When $\eta\in\mathbb{Z}^n$ and $\eta\neq 0$, by  taking the Fourier series in  $d_t v=h$, we obtain
	\begin{equation}\label{equation prop Lz}
		i\sum_{j=1}^{n}\eta_jdt_j\wedge \widehat{v}(\eta,\xi)=\widehat{h}(\eta,\xi).
	\end{equation}
	
	Additionally, since  $\widehat{h}(t,\xi)$ is exact, it follows that $d_t\widehat{h}(t,\xi)=0$. Thus, 
	
	$$
	i\sum_{j=1}^{n}\eta_jdt_j\wedge\widehat{h}(\eta,\xi)=0.
	$$
	
	By Lemma \ref{lemma solution constant coef system},  a solution $ \widehat{v}(\eta,\xi)$ to equation \eqref{equation prop Lz} is given by \eqref{solution lemma}. 	Therefore, the corresponding $v$ is a smooth $p$-form and it is a solution of $d_tv=h$,  which proves the previous statement. 
	
	\bigskip

	Finally, consider the sequence  $\widehat{u}(t,\xi)=\widehat{v}(t,\xi)e^{-i\psi_\xi(t)}$, $\xi\in\mathcal{Z}$, where $\widehat{v}(t,\xi)$ are the Fourier coefficients of the solution obtained previously. Defining the  $p$-form
	$$
	u(t,x)=\sum_{\xi\in\mathcal{Z}} \widehat{u}(t,\xi)e^{i\xi\cdot x}\in C^{\infty}(\mathbb{T}^{n+N};\wedge^{p,0}), 
	$$
	we obtain, for each $\xi\in\mathcal{Z}$,  the following  equation 
	\begin{eqnarray*}
		d_t \widehat{u}(t,\xi)+ic_{\xi0} \wedge \widehat{u}(t,\xi)= e^{-i\psi_\xi(t)}d_t \widehat{v}(t,\xi)= \widehat{f}(t,\xi).
	\end{eqnarray*}
	
	It follows that $\mathbb{L}_{\mathcal{Z}}^pu=f$ and, consequently, $\mathbb{L}_{\mathcal{Z}}^p$ is globally solvable.

\end{proof}

Now, we are ready to state the main result in this section.

\begin{theorem}\label{constant coefficients thm}
	For each $p=0,1,\ldots,n-1$, the operator $\mathbb{L}^p$ is globally solva\-ble if and only if  there exist positive constants  $C$ and $\lambda$ such that 
	\begin{equation}\label{main condition solvability}
		\|\widehat{\mathbb{L}}(\eta,\xi)\|\geq  C|(\eta,\xi)|^{-\lambda},
	\end{equation}
	for all $(\eta,\xi)\in\mathbb{Z}^{n+N}$, with $\widehat{\mathbb{L}}(\eta,\xi)\neq 0$.
\end{theorem}

\begin{proof}
	Assume first that \eqref{main condition solvability} holds. Given 
	$f\in \mathbb{E}^p$ we must \\
	 exhibit $u\in\mathcal{D}'(\mathbb{T}^{n+N};\wedge^{p,0})$ such that 
	$\mathbb{L}^pu=f$, or equivalently, 
	\begin{equation}\label{equation fourier coefficents}
		\widehat{\mathbb{L}}(\eta,\xi)\wedge \widehat{u}(\eta,\xi)=\widehat{f}(\eta,\xi), 
	\end{equation}
	for all $(\eta,\xi)\in\mathbb{Z}^{n+N}$.\\
	
	Since   $\mathbb{L}^{p+1}  f=0$,  we obtain  for each  $(\eta,\xi)\in\mathbb{Z}^{n+N}$  the equation
	\begin{equation}\label{equation fourier homogeneous}
		\widehat{\mathbb{L}}(\eta,\xi)\wedge \widehat{f}(\eta,\xi)=0. 
	\end{equation}
	
	Note that, if $(\eta,\xi)\in\mathbb{Z}^{n+N}$ is such that $\widehat{\mathbb{L}}(\eta,\xi) = 0$, then the constant 1-form $c_{\xi0}$ is integral. In this case,  we take $f_0 \in \mathbb{E}_{0}^p$ defined by 
	$$
	\widehat{f_0}(t,\xi) = 
	\begin{cases}
	\widehat{f}(t,\xi), & \xi \in \mathcal{Z},\\
	0, & 	 \xi \notin \mathcal{Z}.
	\end{cases}
%	\left\{
%	\begin{array}{l}
%	\widehat{f}(t,\xi), \ \xi \in \mathcal{Z}, \\
%	0, \ \xi \notin \mathcal{Z}.
%	\end{array}
%	\right.
	$$	
	Then, in view of Proposition \ref{proposition integral form}, we may define $\widehat{u}(\eta,\xi)$ from $\mathbb{L}_{\mathcal{Z}}^pu=f_0$.

	If  $(\eta,\xi)\in\mathbb{Z}^{n+N}$ is such that  $\widehat{\mathbb{L}}(\eta,\xi)\neq 0$, since $\widehat{f}(\eta, \xi)$ satisfies \eqref{equation fourier homogeneous}, by Lemma \ref{lemma solution constant coef system}, we obtain a solution to  equation  \eqref{equation fourier coefficents} given by
	\begin{align*}
		\widehat{u}(\eta,\xi) &= \frac{1}{i} \underset{\mu \in K}{\sum_{|K|=p+1}} (-1)^{\text{sign}(\mu ,k_1,\ldots,\widehat{\mu },\ldots,k_{p+1})} \cdot \frac{1}{\eta_\mu  + p_\mu (\xi)} \widehat{f_K}(\eta,\xi) dt_{K-\{\mu \}},
	\end{align*}
	where $\eta_\mu  + p_{\mu }(\xi)$ is such that $	|\eta_\mu  + p_{\mu }(\xi)|=\max_{j=1,\dots,n}|\eta_j + p_{j}(\xi)|. 
	$
	
	Since each $f_K$ is a smooth function on 
	$\mathbb{T}^{n+N}$ and  we are assuming the validity of \eqref{main condition solvability}, then 
	for each $\lambda>0$ there exists $C_2>0$ such that
	{$$
		\|	\widehat{u}(\eta,\xi) \|\leq C_1 \frac{1}{|\eta_\mu  + p_\mu (\xi)|} \|\widehat{f}(\eta,\xi)\|\leq C_2(1+|(\eta,\xi)|)^{-\lambda},
		$$}
	for all $(\eta,\xi)\in\mathbb{Z}^{n+N}$,  such that  $\widehat{\mathbb{L}}(\eta,\xi)\neq 0$.

	We thus conclude that   
	$$
	u=\sum_{(\eta,\xi)\in\mathbb{Z}^{n+N}} 	\widehat{u}(\eta,\xi) e^{i(\eta \cdot t+\xi\cdot  x)}\in C^\infty(\mathbb
	T^{n+N};\wedge^{{p},0}), 
	$$ 
	obtained by considering the previously described Fourier coefficients, is a solution   to the equation $\mathbb{L}^pu=f$. 
	
	\bigskip

	Now, suppose that \eqref{main condition solvability} does not hold. Thus, there exists a sequence $(\eta_l,\xi_l)\in\mathbb{Z}^{n+N}$, such that $|(\eta_l,\xi_l)|$ is increasing and 
	$$
	0< \|\widehat{\mathbb{L}}(\eta_l,\xi_l)\|< |(\eta_l,\xi_l)|^{-l}, \quad l\in\mathbb{N}.
	$$

	Let $0<2\delta<1$. For  the case $p=0$, we consider the sequence in $C^\infty(\mathbb
	T^{n+N};\wedge^{0,0})\simeq C^\infty(\mathbb
	T^{n+N})$ given by
	$$
	u_l(t,x)=|(\eta_l,\xi_l)|^{\delta l} e^{i(\eta_l\cdot t+\xi_l \cdot x )},\quad l\in\mathbb{N},
	$$
	while for $p=1,2,\ldots,n-1$, we set the following  sequence in $C^\infty(\mathbb
	T^{n+N};\wedge^{p,0})$ 
	$$
	u_l(t,x)=|(\eta_l,\xi_l)|^{\delta l} e^{i(\eta_l\cdot t+\xi_l \cdot x )} dt_2\wedge \cdots \wedge dt_{p+1},\quad l\in\mathbb{N}.
	$$
	
	For each  $j=1,\ldots,n$, we define  the formal series
	\begin{eqnarray*}\label{defi fj}
		f_j(t,x) =  \sum_{l=1}^\infty i |(\eta_l,\xi_l)|^{\delta l}(\eta_{j,l}+p_j(\xi_l))e^{i(\eta_l \cdot t + \xi_l \cdot x)}.
	\end{eqnarray*} 
	
	Therefore, by hypothesis, we obtain 
	\begin{eqnarray*}
		|\widehat{f_j}(\eta_l,\xi_l)| &=& |(\eta_l,\xi_l)|^{\delta l} |i(\eta_{j,l}+p_j(\xi_l))|%\leq E_l|\eta_l,-\xi_l\alpha_{0}| 
		\\\\ &\leq& |(\eta_l,\xi_l)|^{\delta l}\|\widehat{\mathbb{L}}(\eta_l,\xi_l)\|\\\\&<& |(\eta_l,\xi_l)|^{-l/2}, 
	\end{eqnarray*}
	for all $l\in \mathbb{N}$. As a result $f_j\in C^{\infty}(\mathbb{T}^{n+N})$.\\ 
	
	Consider the $(p+1)$-form in $C^\infty (\mathbb{T}^{n+N}, \wedge^{p+1,0})$ given by 
	$$
	f(t,x)= \sum_{j=1}^n  f_j(t,x) dt_j,\;\text{ when }p=0,
	$$
	and
	$$
	f(t,x)= \sum_{j=1}^n  f_j(t,x) dt_j\wedge dt_2\wedge \cdots \wedge dt_{p+1},\;\text{ when }p=1,2,\ldots,n-1.
	$$

	By the previous construction, we obtain that $\mathbb{L}^{p+1}f=0$. 		Additionally,  if $c_{\xi_l0}=\sum_{j=1}^{n} p_j(\xi_l)dt_j$ is integral, we consider  
	$$
	v_l(t)=|(\eta_l,\xi_l)|^{\delta l}e^{i(\psi_{\xi_l}(t)+\eta_l\cdot t)},\quad \text{when} \;\;p=0;
	$$ 
	$$  
	v_l(t)=  |(\eta_l,\xi_l)|^{\delta l}e^{i(\psi_{\xi_l}(t)+\eta_l\cdot t)} dt_2\wedge \cdots \wedge dt_{p+1},\quad \text{when} \;\; p=1,2,\dots,n-1.
	$$ 
	We conclude that $d_t v_l(t)=e^{i\psi_{\xi_l}(t)}\widehat{f}(t,\xi_l)$ and,  consequently, the $(p+1)$-form $e^{i\psi_{\xi_l}(t)}\widehat{f}(t,\xi_l)$ is exact.  Therefore, $f\in\mathbb{E}^{p}$.\\
	
	Suppose that there exists $v\in \mathcal{D}'(\mathbb{T}^{n+N};\wedge^{p,0}
	)$ satisfying $\mathbb{L}^p v = f$. Thus, by replacing
	the formal Fourier series of $v$ and $f$ in the equation  $\mathbb{L}^p v = f$, we obtain 
	\begin{equation}\label{Fourier coef soluion}
		\widehat{\mathbb{L}}(\eta_l,\xi_l)\wedge\widehat{v}(\eta_l,\xi_l)=\widehat{f}(\eta_l,\xi_l),
	\end{equation}
	for all $l\in\mathbb{N}$. A similar argument applied to the equation $\mathbb{L}^{p+1} f=0$ results in  
	$$
	\widehat{\mathbb{L}}(\eta_l,\xi_l)\wedge \widehat{f}(\eta_l,\xi_l)=0,
	$$
	for all $l\in\mathbb{N}$.

	Moreover, since   $\widehat{u_l}(\eta_l,\xi_l)$ is also a solution to equation  \eqref{Fourier coef soluion}, it follows from Lemma \ref{lemma solution constant coef system}   that  
	\begin{eqnarray}\label{Fourier coef solution Lemma}
		\widehat{u}_l(\eta_l,\xi_l) =   \widehat{v}(\eta_l,\xi_l) + \widehat{\mathbb{L}}(\eta_l,\xi_l)\wedge w_l,
	\end{eqnarray}
	for all $l\in\mathbb{N}$, where $w_l\in \wedge^{p-1}\mathbb{C}^n$.

	Thus, if we consider  the first  level of the complex ($p=0$), we have   $w_l\in \wedge^{-1}\mathbb{C}^n=\{0\}$, then
	$$
	\widehat{v}(\eta_l,\xi_l) = \widehat{u}_l(\eta_l,\xi_l)=|(\eta_l,\xi_l)|^{\delta l},
	$$ 
	for all $l\in\mathbb{N}$. It follows that $\widehat{v}(\eta_l,\xi_l)$ does
	not correspond to any sequence of Fourier coefficients of a distribution $v\in \mathcal{D}'(\mathbb{T}^{n+N})$, which is a contradiction. We thus have that $\mathbb{L}^0$ is not globally solvable.\\

	Consider now the remaining case $p =1,2,\ldots, n-1$. Since there exists $k\in\{1,\ldots,n\}$, such that $ |\eta_{k,l} + p_k(\xi_l)|=\max_{j=1, \ldots, n} |\eta_{j,l} + p_j(\xi_l)|$, for infinitely many $l\in\mathbb{N}$, by passing to a subsequence and reordering the operators $L_j$ in System \eqref{system constant coefficients}, we may assume that 
	\begin{equation}\label{a10 maximum}
		|\eta_{1,l}+p_1(\xi_l)|=\|\widehat{\mathbb{L}}(\eta_l,\xi_l)\|<  |(\eta_l,\xi_l)|^{-l},
	\end{equation}
	for all $l\in\mathbb{N}$.

	In this case, we consider $\zeta_l = \widehat{u}_l(\eta_l,\xi_l) - \widehat{v}(\eta_l,\xi_l)  $
	and  we rewrite  equation \eqref{Fourier coef solution Lemma} omitting the index $l$ to simplify the notation. Also, in order to obtain a contradiction, we only consider the coefficients of
	terms involving $dt_{1},\ldots,dt_{p+1}$ in the  equation $\widehat{\mathbb{L}}(\eta,\xi)\wedge w = \zeta$.  Precisely, we consider the following equality
	\begin{equation}\label{igualdade-necessidade}
		\sum_{j=1}^{p+1}i(\eta_j+p_j(\xi))dt_j \wedge \sum_{|J|=p-1} w_J dt_J= \sum_{|K|=p} \zeta_K dt_K,
	\end{equation}
	where $J= (j_1,\ldots,j_{p-1})$, with $1\leq j_1<j_2<\cdots<j_{p-1} \leq p+1$ and  $K= (k_1,\ldots,k_p)$, with  $1\leq k_1<k_2<\cdots<k_p \leq p+1$.
	
	We need to introduce a little more notation.
	\begin{align*}
		dt^{(j)} &= dt_1\wedge \cdots \wedge \widehat{dt_j} \wedge \cdots \wedge dt_{p+1} ,\\
		\zeta^{(j)} &= \zeta_K, \ \textrm{with} \ K = (1,2,\ldots,\widehat{j},\ldots,p+1),
	\end{align*}
	and, if $j<k$ we write %
	\begin{align*}
		dt^{(j,k)} &= dt_1 \wedge \cdots \wedge \widehat{dt_j} \wedge \cdots \wedge \widehat{dt_k} \wedge \cdots \wedge dt_{p+1}, \\
		w^{(j,k)} &= w_J, \ \textrm{with} \ J=(1,\ldots,\widehat{j},\ldots,\widehat{k},\ldots,p+1).
	\end{align*}
	
	Thus,  from \eqref{igualdade-necessidade}, we obtain the following 
	\begin{align*}
		\sum_{j=1}^{p+1} \zeta^{(j)} dt^{(j)} &= \sum_{k=1}^{p+1} i(\eta_k+p_k(\xi)) dt_k \wedge \sum_{1\leq \mu < \nu \leq p+1} w^{(\mu , \nu)} dt^{(\mu,\nu)} \\
		&= \sum_{k=1}^{p+1} \sum_{1\leq \mu < \nu \leq p+1} i(\eta_k+p_k(\xi)) w^{(\mu , \nu)} dt_k\wedge dt^{(\mu,\nu)} \\
		&= \sum_{k=1}^{p+1} \sum_{1\leq \mu < k} i(\eta_k+p_k(\xi)) w^{(\mu , k)} dt_k \wedge dt^{(\mu,k)} \\&+ \sum_{k=1}^{p+1} \sum_{k < \nu \leq p+1} i(\eta_k+p_k(\xi)) w^{(k , \nu)} dt_k \wedge dt^{(k,\nu)}\\
		& =\sum_{k=1}^{p+1} \sum_{1\leq \mu < k} i(\eta_k+p_k(\xi)) w^{(\mu , k)} (-1)^k dt^{(\mu)} \\&+ \sum_{k=1}^{p+1} \sum_{k < \mu \leq p+1} i(\eta_k+p_k(\xi)) w^{(k , \mu)} (-1)^{k-1} dt^{(\mu)}. 
	\end{align*}
	
	Note that for  $j=1$ we obtain
	$$ \zeta^{(1)}= \sum_{1< k \leq p+1} i(\eta_k+p_k(\xi)) w^{(1 , k)} (-1)^k.
	$$
	
	For $j=2$ it follows that
	\begin{align*}
		\zeta^{(2)} &= i(\eta_1+p_1(\xi)) w^{(1 , 2)} (-1)^0  + \sum_{2< k\leq p+1} i(\eta_k+p_k(\xi)) w^{(2 , k)} (-1)^k\\ &= \sum_{1\leq k <2} i(\eta_k+p_k(\xi)) w^{(k , 2)} (-1)^{k-1} + \sum_{2 < k\leq p+1} i(\eta_k+p_k(\xi)) w^{(2 , k)} (-1)^k .
	\end{align*}
	
	In general, for each $j=1,\ldots,p+1$, we conclude that
	\begin{equation*}
		\zeta^{(j)}  =  \sum_{1\leq k <j} i(\eta_k+p_k(\xi)) w^{(k , j)} (-1)^{k-1} + \sum_{j < k\leq p+1} i(\eta_k+p_k(\xi)) w^{(j , k)} (-1)^k .
	\end{equation*}

	Multiplying the above equality by $(-1)^{j+1} i(\eta_j+p_j(\xi))$ and adding $j$ from 1 to $p+1$, we obtain
	\begin{equation}\label{soma-dos-Aj}
		\sum_{j=1}^{p+1} (-1)^{j+1} i(\eta_j+p_j(\xi))  \zeta^{(j)}  = 0. 
	\end{equation}
	
	Now, returning with index $l$, note that when $K=(2,\ldots,p+1)$ we have  
	$$
	\widehat{v}_{K}(\eta_l,\xi_l) =   |(\eta_l,\xi_l)|^{\delta l}- \zeta_l^{(1)}, 
	$$  and when $K=(1,\ldots,\widehat{j},\ldots,p+1)$, with  $j\neq 1$, we have 
	$$
	\widehat{v}_{K}(\eta_l,\xi_l) = - \zeta_l^{(j)}.
	$$ 
	
	Therefore, since the Fourier coefficients  of  
	$v\in \mathcal{D}'(\mathbb{T}^{n+N};\wedge^{p,0})$ have a slow growth, consequently, there exist positive constants $C_j$ and $\lambda_j$ such that 
	$$
	\big||(\eta_l,\xi_l)|^{\delta l}- \zeta_l^{(1)}\big|   \leq C_1 { |(\eta_l,\xi_l)|^{\lambda_1}},
	$$
	for all $l\in \mathbb{N}$, and 
	$$
	|\zeta_l^{(j)}|\leq C_j{ |(\eta_l,\xi_l)|^{\lambda_j}},\quad j=2,\ldots,p+1,
	$$
	for all $l\in \mathbb{N}$.

	Since $\eta_{1l}+p_1(\xi_l)\neq 0$, it follows from  \eqref{soma-dos-Aj} that 
	$$
	\zeta_l^{(1)} = \sum_{j=2}^{p+1} (-1)^j \frac{\eta_{j,l}+p_j(\xi_l)}{\eta_{1,l}+p_1(\xi_l)}\zeta_l^{(j)}.
	$$
	Therefore, we have the following equation
	\begin{eqnarray*}
		|(\eta_l,\xi_l)|^{\delta l}&=& \left(|(\eta_l,\xi_l)|^{\delta l}-\zeta_l^{(1)}\right) +\zeta_l^{(1)} \\
		&=& \left(|(\eta_l,\xi_l)|^{\delta l}-\zeta_l^{(1)}\right) + \sum_{j=2}^{p+1} (-1)^{j} \frac{\eta_{j,l}+p_j(\xi_l)}{\eta_{1,l}+p_1(\xi_l)}\zeta_l^{(j)}.
	\end{eqnarray*}

	We recall that, by \eqref{a10 maximum}, we have  $|\eta_{j,l}+p_j(\xi_l)|\leq |\eta_{1,l}+p_1(\xi_l)|$, for all $j=1,\ldots,p+1$, and thus  we obtain the following inequality 
	\begin{eqnarray*}
		|(\eta_l,\xi_l)|^{\delta l}&\leq &\big||(\eta_l,\xi_l)|^{\delta l} - \zeta_l^{(1)}\big| + \sum_{j=2}^{p+1} |\zeta_l^{(j)}| \\
		&\leq& C_1 { |(\eta_l,\xi_l)|^{\lambda_1}}+ \sum_{j=2}^{p+1} C_j{ |(\eta_l,\xi_l)|^{\lambda_j}}\\ &\leq& C (p+1) |(\eta_l,\xi_l)|^{\lambda},
	\end{eqnarray*}
	where  $C = \underset{1\leq j\leq p+1}{\max} C_j$ and $\lambda=\underset{1\leq j\leq p+1}{\max} \lambda_j$.\\
	
	The previous inequality implies that the sequence  $ |(\eta_l,\xi_l)|^{\delta l-\lambda}$, $l\in\mathbb{N}$,  is  bounded,  which is a  contradiction.

\end{proof}

Let $\alpha=\sum_{j=1}^{n}\alpha_jdt_j$ be a real, closed and  constant 1-form on $\mathbb{T}^n$. If $\alpha$ is not rational, then Definition 2.1, together with Proposition 2.2 in \cite{BerCorMal93}, imply that    $\alpha$ is a Liouville form if and only if for each $\lambda>0$, there exist a constant  $C>0$ and a sequence $(\eta_l,\xi_l)$ in $\mathbb{Z}^{n+N}$, $|\xi_l|\geq 2$,  such 
$$
\max_{j=1, \ldots, n} |\eta_{j,l}+\alpha_j\xi_l|\leq C |(\eta_l,\xi_l)|^{-l}, \quad l\in\mathbb{N}. 
$$

As a consequence of Theorem \ref{constant coefficients thm}, we recover the main result in \cite{BerPet99jde}. Precisely,  for each $j=1,\ldots,n$,   consider the operador  $P_j(D_x)=D_x$ in System \eqref{main system} and  $c_j$  a real valued function; then, we have the following result.

\begin{theorem}\label{differential_complex}
	For each $p=0,1,\ldots,n-1$, the operator $$
	\mathbb{L}^p=d_t+a_0\wedge \dfrac{\partial}{\partial x},
	$$
	defined on $\mathbb{T}_t^{n}\times \mathbb{T}_x^{1} $, is globally solvable if and only if the real and closed 1-form $a_0=\sum_{j=1}^{n}a_{j0}dt_j$ is either rational or non-Liouville. 
\end{theorem}

\begin{proof}
	In fact, if $a_0$ is  not rational, then $\widehat{\mathbb{L}}(\eta,\xi)=\sum_{j=1}^{n}i(\eta_j+a_{j0}\xi)dt_j\neq 0$, for all $(\eta,\xi)\neq 0$. Thus, Condition \eqref{main condition solvability} can be written as 
	follows: there exist positive constants  $\lambda$ and  $C$ such that 
	$$
	\|\widehat{\mathbb{L}}(\eta,\,\xi)\|=\max_{j=1, \ldots, n} |\eta_j + a_{j0}\xi| \geq  C|(\eta,\xi)|^{-\lambda}, 
	$$
	for all $(\eta,\xi)\neq (0,0)$, or equivalently, $a_0$ is non-Liouville.

	\bigskip
	
	On the other hand, when  $a_0$ is rational, there exists a minimal natural number $q$ such that $qa_0$ is integral. Note that, if $q=1$, then $c_{\xi0}=\xi a_0 $
	is integral for all $\xi\in\mathbb{Z}$ and the result follows from Proposition \ref{proposition integral form}. 
	
	In the case $q>1$, we claim that there exists $C>0$ such that 
	$$
	\max_{j=1, \ldots, n} |\eta_j + a_{j0}\xi| \geq C,$$
	for all $(\eta,\xi)\in\mathbb{Z}^{n+1}$, when   $\widehat{\mathbb{L}}(\eta,\xi)\neq0$. 
	
	In fact, if $\widehat{\mathbb{L}}(\eta,\xi)\neq0$ and 	$\xi\in (q\mathbb{Z}),$ then 	$
	\max_{j=1, \ldots, n} |\eta_j + a_{j0}\xi| \geq 1$. Moreover, if $\xi\notin(q\mathbb{Z})$, then there exist $d\in
	\mathbb{Z}$ and $r\in \{1,\ldots,q-1\}$ such that $\xi=dq+r$. Thus, we have
	\begin{eqnarray*}
		\max_{j=1, \ldots, n} |\eta_j + a_{j0}\xi|&=& r\max_{j=1, \ldots, n} \left|\dfrac{\eta_j+dqa_{j0}}{r} + a_{j0}\right|\geq dist(\alpha_0, r^{-1}\mathbb{Z}^n)\geq C_0,
	\end{eqnarray*}
	where $\alpha_0=(a_{10},\ldots,a_{n0})$ and  $C_0=\min_{r\in\{1,\ldots,q-1\}} dist(\alpha_0, r^{-1}
	\mathbb{Z}^n)>0$.	Thus, the conclusion of the previous statement follows by taking $C=\min\{1,C_0\}$. Therefore, condition \eqref{main condition solvability} holds and the proof is complete.

\end{proof}

\begin{remark}\label{remark coef constants}
	As a consequence of Theorem \ref{constant coefficients thm},  if there exists $k\in\{1,\ldots,n\}$ such that $L_k=D_{t_k}+P_k(D_x)$, defined on $\mathbb{T}_t^1\times
	\mathbb{T}_x^N$, is globally solvable, then the operator  $\mathbb{L}^p$ given in \eqref{operatorconstant coefficients} is also globally solvable. Indeed, if $\mathbb{L}^p$  were not globally solvable, there would be a sequence $(\eta_l,\xi_l)\in\mathbb{Z}^{n+N}$, such that
	$$
	|\eta_{l}^{(k)}+p_k(\xi_l)|\leq \|\widehat{\mathbb{L}}(\eta_l,\xi_l)\|< |(\eta_l,\xi_l)|^{-l}, \quad l\in\mathbb{N}, 
	$$
	concluding that $L_k$ is also not globally solvable, which ensures the statement.

	However, the converse of this result is not true. In fact, in  Example 4.9 in \cite{Ber99}, two Liouville numbers $\alpha_1 $ and $\alpha_2$  are presented such that $(\alpha_1,\alpha_2)$ is a non-Liouville vector. Consider the following system defined on $\mathbb{T}^3$
	$$
	L_j=D_{t_j}+\alpha_j D_x, \quad j=1,2.
	$$
	In this case, although neither $L_1$ nor $L_2$ is  globally solvable in $\mathbb{T}^2$, we have  that  $a_0=\alpha_1dt_1+\alpha_2dt_2$ is a real, closed and non-Liouville form. Consequently,  by  Theorem \ref{differential_complex}, the associated operator $\mathbb{L}^p$ is globally solvable.
\end{remark}

%====================================================
%====================================================
\subsection{Homogeneous symbols \label{homo-symbols}}
%====================================================
%====================================================

In this section we investigate the global solvability  of  $\mathbb{L}^p$ associated with the following special class of systems 
\begin{equation}\label{homogeneous system}
	L_j=D_{t_j}+P_j(D_x),\quad j=1,\ldots,n,
\end{equation}
on $\mathbb{T}_t^n\times \mathbb{T}_x^1$, where for each $j=1,\ldots,n$,  the  symbol of the  pseudo-differential operator $P_j(D_x)$  is  given by
$$
p_j(\xi)=c_j|\xi|^{\kappa}, \quad c_j\in\mathbb{C},\; \xi\in \mathbb{Z}_*,
$$
where $\kappa>0$ is a rational number. In particular, $p_j$ is a homogeneous  symbol of order $\kappa$ (see 
\cite{AviGonKirMed19}). 

In this case, by
a slight modification in Condition \eqref{main condition solvability} in Theorem \ref{constant coefficients thm}, we have that $\mathbb{L}^p$ is globally solvable
if and only if there exist positive constants $C$ and $\lambda$ such that the following  inequality holds 
\begin{equation}\label{condition homogeneous symbol}
	\max_{j=1,\ldots,n}\left|\dfrac{\eta_j}{\xi^\kappa}+c_j\right|\geq C (|\eta|+\xi)^{-\lambda}, 
\end{equation}
for all $(\eta,\xi)\in\mathbb{Z}^{n}\times \mathbb{N}$, such that $\widehat{\mathbb{L}}(\eta,\xi)=i\sum_{j=1}^{n}(\eta_j+c_j\xi^\kappa)dt_j \neq 0$. \bigskip

We will make use of the notations: $c_j=\alpha_j+i\beta_j$, with $\alpha_j,\beta_j\in\mathbb{R}$ and   $\alpha=(\alpha_1,\ldots,\alpha_n),\beta=(\beta_1,\ldots,\beta_n)\in\mathbb{R}^n$. Also, taking inspiration from \cite{AviMed21}, we  recall the following definition about simultaneous Diophantine approximations $(SDA)$ of real numbers by rationals with common  denominators.
\begin{definition}\label{def-SDA}	
	Let $\alpha=(\alpha_1,\ldots,\alpha_n)\in\mathbb{R}^n\setminus \mathbb{Q}^n$ be a vector and $\mu\in\mathbb{N}$. We say that $\alpha$ satisfies the condition $(SDA)_{\mu}$ if there exist a positive  constant $C$ and a sequence $(p_\ell,q_\ell)=(p_{1\ell},\ldots,p_{n\ell},q_\ell)\in\mathbb{Z}^n\times \mathbb{N}$ such that 
	\begin{equation}
		\max_{j=1,\dots,n}\left|\alpha_{j}^\mu-\dfrac{(p_{j\ell})^\mu}{q_\ell}\right|<  C q_{\ell}^{-\ell}, \quad \ell \in\mathbb{N}.
	\end{equation}
	In particular, if $\mu=1$, we say that $\alpha$ is a Liouville vector (see \cite{BerMedZan12}).
\end{definition}

When $\mu>1$, there are vectors $\alpha\in\mathbb{R}^n\setminus \mathbb{Q}^n$ satisfying the $(SDA)_\mu$ condition, which are non-Liouville vectors (see Example  2 in \cite{AviMed21}).

\begin{theorem}\label{homogeneous operators}
	Let $\kappa={\rho}/{\mu}$ be a positive rational number such that  $\rho,\mu\in\mathbb{N}$ and $\gcd(\rho,\mu)=1$. 
	For each $p=0,1,\ldots,n-1$, the operator $\mathbb{L}^p$ associated with
	System \eqref{homogeneous system}	is globally solvable if and only if  $\alpha\in\mathbb{R}^n$ is either rational or it does not satisfy $(SDA)_{\mu}$, whenever $\beta=0$.	
\end{theorem}

\begin{proof} 
	When $\beta\neq 0$  the validity of \eqref{condition homogeneous symbol} follows immediately and, consequently, $\mathbb{L}^p$ is globally solvable. Thus, from now on, we assume that $\beta=0$.

	Consider first $\alpha\in\mathbb{R}^n\setminus\mathbb{Q}^n$. If $\mathbb{L}^p$ is not globally solvable,
	it follows from \eqref{condition homogeneous symbol} that there exists a sequence $(\eta_{\ell},\xi_{\ell})=(\eta_{1\ell},\ldots,\eta_{n\ell},\xi_{\ell})\in\mathbb{Z}^n\times\mathbb{N}$ such that
	\begin{equation}\label{equation thm homog}\max_{j=1,\ldots,n}\left|\dfrac{\eta_{j\ell}}{\xi_{\ell}^{\kappa}}-\alpha_j\right|<(|\eta_{\ell}|+\xi_{\ell})^{-\rho\ell},\;\;\ell\in\mathbb{N}.
	\end{equation}
	
	For each $j=1,\ldots,n$, we  write 
	\begin{equation}\label{nice equality}
		\left|\dfrac{(\eta_{j\ell})^\mu}{\xi_\ell^{\rho}}-\alpha_j^\mu\right|= \left|\dfrac{\eta_{j\ell}}{\xi_\ell^{\kappa}}-\alpha_j\right|\cdot \left|\sum_{r=1}^{\mu}\left(\dfrac{\eta_{j\ell}}{\xi_\ell^{\kappa}}\right)^{\mu-r}(-\alpha_j)^{r-1}\right|,\quad \ell\in\mathbb{N}.
	\end{equation}

	By \eqref{equation thm homog}, the sequence $\sum_{r=1}^{\mu}\left|({\eta_{j\ell}}/{\xi_{\ell}^{\kappa}})^{\mu-r}\alpha_j^{r-1}\right|$   converges to $\mu|\alpha_j|^{\mu-1}$ as $\ell$ goes to infinity. Thus, there exists $C>0$ such that    
	\begin{eqnarray*}
		\max_{j=1,\ldots,n}\left|\dfrac{(\eta_{j\ell})^{\mu}}{\xi_{\ell}^{\rho}}-\alpha_j^\mu\right|
		\leq C (|\eta_{\ell}|+\xi_{\ell})^{-\rho\ell}\leq C (\xi_{\ell}^{\rho})^{-\ell},\quad \ell\in \mathbb{N}.
	\end{eqnarray*}
	Therefore,  $\alpha$ satisfies the  $(SDA)_\mu$ condition.
	
	\bigskip
	
	Supposing that $\alpha\in\mathbb{R}^n\setminus\mathbb{Q}^n$ satisfies  $(SDA)_\mu$, there exist a positive constant $C_0$ and a sequence $(p_\ell,q_\ell)=(p_{1\ell},\ldots,p_{n\ell},q_\ell)\in\mathbb{Z}^n\times \mathbb{N}$  satisfying the following inequality
	$$
	\max_{j=1,\ldots, n}	\left|\alpha_j^{\mu}-\dfrac{(p_{j\ell})^ {\mu}}{q_\ell}\right|\leq C_0 (|p_\ell|+q_\ell)^{-\ell},\quad \; \ell\in\mathbb{N}.
	$$	
	In particular, if there exists $\alpha_j=0$, then $p_{j\ell}=0$, for large enough $\ell$.
	
	\bigskip
	
	Additionally, since gcd$(\rho,\mu)=1$, there exist $\widetilde{\rho}, \widetilde{\mu}\in\mathbb{N}$ such that $\widetilde{\rho}\rho-\widetilde{\mu}\mu=1$. Thus, for each  $j=1,\ldots,n$, such that $\alpha_j\neq0$, we obtain:
	\begin{eqnarray*}
		C_0(|p_\ell|+q_\ell)^{-\ell}&\geq& 	\left|\dfrac{(p_{j\ell})^{\mu}}{q_\ell}-\alpha_j^{\mu}\right | \\&=& \left|\dfrac{(p_{j\ell} q_{\ell}^{\widetilde{\mu}})^{\mu}}{q_\ell^{\rho\widetilde{\rho}}}-\alpha_j^{\mu}\right|  
		\\&=& \left|\dfrac{p_{j\ell} q_{\ell}^{\widetilde{\mu}}}{(q_{\ell}^{\widetilde{\rho}})^{\kappa}}-\alpha_j\right| \cdot  \left|\sum_{r=1}^{\mu} \left(\dfrac{p_{j\ell} q_{\ell}^{\widetilde{\mu}}}{(q_{\ell}^{\widetilde{\rho}})^{\kappa}}\right)^{\mu-r}(-\alpha_j)^{r-1}\right| \\
		& = &\left|\dfrac{p_{j\ell} q_{\ell}^{\widetilde{\mu}}}{(q_{\ell}^{\widetilde{\rho}})^{\kappa}}-\alpha_j\right| \cdot  \left|\sum_{r=1}^{\mu} \left(\dfrac{(p_{j\ell})^{\mu} }{q_{\ell}}\right)^{({\mu-r})/\mu}\alpha_j^{r-1}\right|\\\
		&\geq & C_j \left|\dfrac{p_{j\ell} q_{\ell}^{\widetilde{\mu}}}{(q_{\ell}^{\widetilde{\rho}})^{\kappa}}-\alpha_j\right|, 
	\end{eqnarray*}
	for large enough $\ell$, where $C_j>0$ can be obtained  since   
	$$
	\lim_{\ell\to \infty}	\left[\sum_{r=1}^{\mu} \left(\dfrac{(p_{j\ell})^{\mu} }{q_{\ell}}\right)^{({\mu-r})/\mu}\alpha_j^{r-1}\right]=\mu\alpha_j^{\mu-1}\neq0.
	$$

	Finally, if $\mathbb{L}^p$ were globally solvable, we would have,  by the previous arguments, together with \eqref{condition homogeneous symbol}, the following inequality
	\begin{eqnarray*}
		{C}_1 C_0 (|p_\ell|+q_\ell)^{-\ell}
		&\geq& \max_{\stackrel{j=1,\ldots,n}{\alpha_j\neq 0}} \left|\dfrac{p_{j\ell} q_{\ell}^{\widetilde{\mu}}}{(q_{\ell}^{\widetilde{\rho}})^{\kappa}}-\alpha_j\right|\\ &=&\max_{j=1,\ldots, n}\left|\dfrac{p_{j\ell} q_{\ell}^{\widetilde{\mu}}}{(q_{\ell}^{\widetilde{\rho}})^{\kappa}}-\alpha_j\right|\\&\geq & C(q_{\ell}^{\widetilde{\mu}}|p_\ell|+q_{\ell}^{\widetilde{\rho}})^{-\lambda},
	\end{eqnarray*}
	where  $\displaystyle{C}_1=\max_{\stackrel{j=1,\ldots,n}{\alpha_j\neq 0}}C_j^{-1}>0$ and $\ell$ is large enough.

	Hence, the sequence  ${(|p_\ell|+q_\ell)^{\ell}}{(q_{\ell}^{\widetilde{\mu}}|p_\ell|+q_{\ell}^{\widetilde{\rho}})^{-\lambda}}$ is bounded, which  is a contradiction. Therefore, $\mathbb{L}^p$ is not globally solvable.
	
	\bigskip
	
	Suppose now that $\alpha\in\mathbb{Q}^n$.  Then, there exists a minimal natural number $q_0$, such that $q_0\alpha\in\mathbb{Z}^n$. In this case, if $\alpha_\mu\doteq (\alpha_1^{\mu},\ldots,\alpha_n^{\mu})\in\mathbb{Q}^n $, then  $q=q_0^\mu$ is the minimal natural number such that $q\alpha_\mu\in\mathbb{Z}^n.$

	Consider the sets  
	$$
	\mathcal{Z}=\left\{\xi\in\mathbb{Z};\; c_{\xi0}=\sum_{j=1}^{n}\alpha_j |\xi|^{\rho/\mu}dt_j \text{ is integral} \right\}
	$$
	and $\mathcal{W}=\mathbb{Z}\setminus\mathcal{Z}$.

	Inspired by Section 4 in \cite{BerCorPet96}, since $\mathcal{Z}\cup \mathcal{W} = \mathbb{Z}$ and  $\mathcal{Z}\cap \mathcal{W}= \emptyset$, the operator  $\mathbb{L}^p$ is globally solvable if and only if both operators 
	$$\mathbb{L}_{\mathcal{Z}}^p \doteq \mathbb{L}^p \vert_{\mathcal{D}_{\mathcal{Z}}'(\mathbb{T}^{n+1};\wedge^{p,0})}\text{ \;and\; }\mathbb{L}_{\mathcal{W}}^p \doteq \mathbb{L}^p \vert_{\mathcal{D}_{\mathcal{W}}'(\mathbb{T}^{n+1};\wedge^{p,0})}
	$$ are globally solvable.

	In Proposition \ref{proposition integral form}, it was proved  that $\mathbb{L}^p_{\mathcal{Z}}$ is globally solvable. If $q=1$, then the proof is complete. When $q>1$, we claim that 
	$\mathbb{L}^p_{\mathcal{W}}$ is also globally solvable. Indeed,  consider the partition $\mathcal{W}=(\mathcal{W}\cap(q\mathbb{Z}))\cup (\mathcal{W}\setminus
	(q\mathbb{Z}))$. 
	
	If $(\eta,\xi)\in\mathbb{Z}^n\times\mathbb{Z}$, where  $\xi\in\mathcal{W}\cap(q\mathbb{Z})$, then $\xi=qm$, for some integer $m\neq 0$, and  
	\begin{eqnarray*}
		|\xi|\max_{j=1,\ldots,n}\left|\alpha_j^{\mu}-\dfrac{\eta_j^{\mu}}{\xi}\right|= \max_{j=1,\ldots,n}\left|\xi\alpha_j^{\mu}-{\eta_j^{\mu}}\right|
		= \max_{j=1,\ldots,n}\left|mq_0^{\mu}\alpha_j^{\mu}-{\eta_j^{\mu}}\right|\geq1, 
	\end{eqnarray*}
	where  the previous inequality occurs because  $(mq_0^{\mu}\alpha_j^{\mu}-{\eta_j^{\mu}})\in\mathbb{Z}_*$, for some $j$. Otherwise, $(mq_0^{\mu}\alpha_j^{\mu}-{\eta_j^{\mu}})=0$ for all $j$, which implies that $m=m_0^\mu$,  with $m_0\in\mathbb{Z}$. This would imply  $\xi=(q_0m_0)^\mu\in\mathcal{Z}$, which is a contradiction. It follows from the previous inequality that 
	\begin{equation}\label{eq rationa case 1}
		\max_{j=1,\ldots,n}\left|\alpha_j^{\mu}-\dfrac{\eta_j^{\mu}}{\xi}\right|\geq \dfrac{1}{|\xi|}.
	\end{equation}
	
	If $(\eta,\xi)\in\mathbb{Z}^n\times\mathbb{Z}$, where   $\xi\in\mathcal{W}\setminus
	(q\mathbb{Z})$, then there exists $d\in\mathbb{Z}$ such that $\xi=dq+r$, where $1\leq r\leq q-1$ ($r\in\mathbb{Z}$).  Thus, we obtain
	\begin{eqnarray*}
		|\xi|\max_{j=1,\ldots,n}\left|\alpha_j^{\mu}-\dfrac{\eta_j^{\mu}}{\xi}\right|&=& \max_{j=1,\ldots,n}\left|\xi\alpha_j^{\mu}-{\eta_j^{\mu}}\right|
		\\\\&=& \max_{j=1,\ldots,n}\left|r\alpha_j^{\mu}+\alpha_j^{\mu}dq-\eta_j^{\mu}\right|\\\\
		&\geq& \max_{j=1,\ldots,n}\left|\alpha_j^{\mu}-\dfrac{(\eta_j^{\mu}-\alpha_j^{\mu}dq)}{r}\right|   \\\\
		&\geq&dist(\alpha_\mu,r^{-1}\mathbb{Z}^n)\geq \mathfrak{D},
	\end{eqnarray*}
	where $\displaystyle \mathfrak{D}=\min_{r=1,\ldots,q-1} dist(\alpha_\mu,r^{-1}\mathbb{Z}^n)$ is a positive constant. It follows that
	\begin{equation}\label{eq rationa case 2}
		\max_{j=1,\ldots,n}\left|\alpha_j^{\mu}-\dfrac{\eta_j^{\mu}}{\xi}\right|\geq \dfrac{\mathfrak{D}}{|\xi|}.
	\end{equation}
	
	Therefore, if  $\widetilde{\mathfrak{D}}=\min\{1,\mathfrak{D}\}$, then, by \eqref{eq rationa case 1} and \eqref{eq rationa case 2}, we obtain 
	\begin{equation*}
		\max_{j=1,\ldots,n}\left|\alpha_j^{\mu}-\dfrac{\eta_j^{\mu}}{\xi}\right|\geq \dfrac{\widetilde{\mathfrak{D}}}{|\xi|},
	\end{equation*}
	for all $(\eta,\xi)\in\mathbb{Z}^{n}\times\mathbb{Z}$, where $\xi\in\mathcal{W}$. Thus, when  $\xi\in \mathcal{W}$,   vector $\alpha$ behaves as if it didn't satisfy the  condition $(SDA)_\mu$, consequently, $\mathbb{L}_{\mathcal{W}}^p$ is globally solvable.
	
\end{proof}

\begin{example} Let $\mu\geq 2$ be an integer. 
	Consider the non-Liouville vector $\alpha=(\sqrt[\mu ]{3/2}\mathscr{L},\sqrt[\mu ]{3\cdot2^{\mu -1}}\mathscr{L})$ in $\mathbb{R}^2$, where $\mathscr{L}=\sum_{k=1}^{\infty}10^{-k!}$ is the Liouville constant (see  Example 2 in \cite{AviMed21}). Thus, $\alpha$ satisfies $(SDA)_\mu$ and, therefore, 
	the operator $\mathbb{L}^p$ associatted with the system 
	$$
	L_j=D_{t_j}+\alpha_{j}(D_x^{2})^{1/2\mu },\quad j=1,2,
	$$
	in $\mathbb{T}_t^2\times\mathbb{T}_x^1$	is not globally solvable. However,  for the same $\alpha$, the operator $\mathbb{L}^p$ associatted with the system 
	$$
	L_j=D_{t_j}+\alpha_{j}(D_x^{2})^{1/2 },\quad j=1,2,
	$$
	in $\mathbb{T}_t^2\times\mathbb{T}_x^1$ is globally solvable.
\end{example}

%====================================
%====================================
\section{Reduction to the normal form}\label{section reduction normal form}
%====================================
%====================================

In this section, we  investigate some connections  between the global sol\-va\-bi\-lity of the operator  $\mathbb{L}^p$, given in  \eqref{main operator}, and the global solvability of its normal form $\mathbb{L}_0^{p}$, which is  associated with  the following  system  with constant coefficients  
$$
L_{j0} = D_t + c_{j0}P_j(D_x), \quad j=1,\ldots, n,
$$
where
$$
c_{j0} = \dfrac{1}{2\pi}\int_{\mathbb{T}^1}c_j(0, \ldots, t_j, \ldots, 0)dt_j.
$$

In this case, we can write the constant coefficient operator $\mathbb{L}_0^p$ as follows
\begin{equation}\label{normal form operator}
	\mathbb{L}^p_0=d_t+c_0(D_x)\wedge,
\end{equation}
where for each   $u=\sum_{|k|=p}u_{K}(t,x)dt_K$, the operator $c_0(D_x)$ is given by 
$$
c_0(D_x)\wedge u=\sum_{j=1}^{n}\sum_{|K|=p}ic_{j0}P_j(D_x)u_Kdt_j\wedge dt_K.
$$

Moreover, the associated  compatibility conditions set, referred as $\mathbb{E}^p_{0}$, is given by the elements
$f\in C^\infty(\mathbb T^{n+N};\wedge^{ {p+1},0})$ such that $\mathbb{L}^{p+1}_0  f=0$ and 
\begin{eqnarray*}\label{cond compat_0}
	\widehat{f}(t,\xi)e^{i\psi_\xi(t)} \;\textrm{is a  $(p+1)$-exact form,}
\end{eqnarray*}
whenever $c_{\xi0} = \sum_{j=1}^{n}c_{j0}p_j(\xi) dt_j$ is integral. The function   $\psi_{\xi}\in{C}^{\infty}(\mathbb{R}^n;\mathbb{R})$ such that $d\psi_\xi=\Pi^*( c_{\xi0})$ is described in Section \ref{comp cond}.\\

This kind of technique,  called reduction to normal form, is widely employed for several classes of linear operators defined on the torus, as seen in \cite{BerPet99jmaa, BerMedZan21,   AriKirMed19, Pet05}. Precisely, the approach consists in obtaining  an automorphism $\Psi$ of the spaces $\mathcal{D}'(\mathbb{T}^{n+N}, \wedge^{p,0})$  and $C^{\infty}(\mathbb{T}^{n+N}, \wedge^{p,0})$ that satisfies the  conjugation formula
\begin{equation}\label{conjugation}
	\Psi\circ \mathbb{L}^p_0 \circ \Psi^{-1} = \mathbb{L}^p,
\end{equation}
where, additionally, $\Psi$ maps $\mathbb{E}^p$ isomorphically onto $\mathbb{E}^p_0$. By assuming the existence of such an automorphism,  if $\mathbb{L}^{p}_0$ is globally solvable, and 
$f \in \mathbb{E}^p$, then $g = \Psi^{-1}f \in \mathbb{E}^p_0$ and there exists $u \in \mathcal{D}'(\mathbb{T}^{n+N}, \wedge^{p,0})$, such that $\mathbb{L}_0^pu = g$. Then, it follows from \eqref{conjugation} that $\mathbb{L}^p(\Psi u) = f$, which implies the global sol\-va\-bi\-lity of $\mathbb{L}^p$. With a similar argument, we may conclude that
the converse is also true.

Now, we state a condition that guarantees the reduction to the normal form of the operator $\mathbb{L}^p$. The isomorphism $\Psi$ that fulfills the conditions described previously will be given in Theorem \ref{reduction normal form thm}.

\begin{definition}\label{definition condition D}
	Let $\alpha=(\alpha_\xi)_{\xi\in\mathbb{Z}^N}$ be  a sequence of complex and    closed 1-forms in $C^{\infty}
	(\mathbb{T}^n;\wedge^{1})$. 
	Consider for each $\xi$  the following decomposition 
	\begin{equation*}\label{exat_reduction}		
		\alpha_\xi = \alpha_{\xi0} + d_t\mathscr{A}_{\xi},
	\end{equation*}
	where $\alpha_{\xi0}\in\wedge^1\mathbb{C}^n$ and   $\mathscr{A}_{\xi} \in C^{\infty}(\mathbb{T}^n)$.

	We say that  $\alpha$ satisfies the condition $\mathscr{D}$, when there exist positive constants $C$ and  $\kappa$ such that
	\begin{equation}\label{condition D}
		\sup_{t\in \mathbb{T}^n}\left\{ 
		\exp \left( \Im \mathscr{A}_{\xi}(t) \right) \right\} \leq 
		C |\xi|^{\kappa}, \ \forall \xi\in\mathbb{Z}^N.
	\end{equation}
\end{definition}

\begin{example}
	Let $c_\xi(t)=\sum_{j=1}^{n}p_j(\xi)c_j(t_j)dt_j$ be a closed 1-form, where each $c_j$ is a smooth function on $\mathbb{T}^1$ and $p_j(\xi)=\mathcal{O}(\log(|\xi|))$, that is,  $p_j(\xi)$ satisfies  \eqref{log-estimate-1}. Then, the sequence  $( c_\xi)_{\xi\in\mathbb{Z}^N}$ satisfies condition $\mathscr{D}$.
\end{example}

The main result in this section reads as follows:
\begin{theorem}\label{reduction normal form thm}
	Assume that the sequence of closed 1-forms $$
	c(t,\xi)=\sum_{j=1}^{n}c_j(t)p_j(\xi)dt_j,\quad\xi\in\mathbb{Z}^N,
	$$ 
	in $C^{\infty}(\mathbb{T}^n;\wedge^1)$ satisfies condition $\mathscr{D}$. Then, for each $p=0,1,\ldots,n-1$, 
	the operator $\mathbb{L}^p$ is globally solvable if and only if the same is true for its normal form $\mathbb{L}^p_0$, given in \eqref{normal form operator}.
\end{theorem}

\begin{proof}
	Consider the decomposition
	$c(t,\xi)=c_{\xi0}+d_t \mathscr{C}_\xi(t)$,
	where $c_{\xi0}\in\wedge^1\mathbb{C}^n$ and  $\mathscr{C}_\xi $ is a sequence of functions in $ C^{\infty}(\mathbb{T}^n)$  that satisfies  \eqref{condition D}, for all $\xi\in\mathbb{Z}^N$.
	
	We claim that, for each $p=0,1,\ldots,n-1$, the following linear operator defines an isomorphism on both spaces $\mathcal{D}'(\mathbb{T}^{n+N}, \wedge^{p,0})$  and $C^{\infty}(\mathbb{T}^{n+N}, \wedge^{p,0})$
	\begin{equation*}
		u=\sum_{\xi\in\mathbb{Z}^N} \widehat{u}(t,\xi) e^{i x \cdot  \xi} \mapsto \Psi u= \sum_{\xi \in \mathbb{Z}^N}  \widehat{u}(t, \xi)
		e^{-i \mathscr{C}_{\xi}(t)e^{i x \cdot  \xi} }.
	\end{equation*}

	Moreover, the inverse is given by	
	\begin{equation*}
		\Psi^{-1} u = \sum_{\xi \in \mathbb{Z}^N} \widehat{u}(t, \xi)
		e^{i \mathscr{C}_{\xi}(t) } e^{i x \cdot  \xi}.
	\end{equation*}
	
	Consider first 
		$$
	u=\sum_{\xi\in \mathbb{Z}^N}
	\sum_{|K|=p}
	\left(
	\widehat{u_K}(t,\xi)e^{i\xi \cdot x}
	\right)dt_K \in C^{\infty}(\mathbb{T}^{n+N}, \wedge^{p,0}),
	$$
	where
	$\widehat{u_K}(t,\xi)$ stands as the Fourier transform of $u_K$.
	
	Since $\widehat{\Psi u}(t,\xi) = \sum_{|K|=p}e^{-i \mathscr{C}_{\xi}(t)}\widehat{u_K}(t,\xi)dt_K$, for each multi-index $K$, with $|K|=p$,
	we must show that for any $\alpha\in\mathbb{Z}_+^n$ and $M\in\mathbb{N}$, there exists  $C> 0$ such that
	\begin{equation}\label{estimative derivatives}
		\left|\partial_{t}^{\alpha} \left( \widehat{u_K}(t,\xi) 	e^{-i \mathscr{C}_{\xi}(t) } \right)\right|\leq C |\xi|^{-M},
	\end{equation}
	for all  $(t,\xi) \in\mathbb{T}_t^n\times  \mathbb{Z}^n$. \\

	Given $\alpha \in \mathbb{Z}_+^n$, we obtain from Leibniz's rule that
	\begin{equation}\label{Leibniz rule}
		\partial_{t}^{\alpha}  \left( \widehat{u_K}(t,\xi) 	e^{-i \mathscr{C}_{\xi}(t) } \right) = \sum_{\beta\leq \alpha}
		\binom{\alpha}{\beta} 
		\partial_{t}^{\beta} \left( e^{-i \mathscr{C}_{\xi}(t)}\right) \, \partial_{t}^{\alpha-\beta}\widehat{u_K}(t,\xi).
	\end{equation}
	Now, if $\beta\in \mathbb{Z}_+^n$ is nonzero, then  by  Fa\`a di Bruno's multivariate  formula (see \cite{ConstantineSavits}), we get
	\begin{equation}\label{Faa Di Bruno Formula}
		\partial_{t}^{\beta} e^{-i \mathscr{C}_{\xi}(t)}  = 
		\sum_{\stackrel{1\leq k \leq |\beta|}{ S(\beta,k)}}\dfrac{\beta!}{k_1! \cdots k_{|\beta|}!}   e^{-i \mathscr{C}_{\xi}(t)}  \prod_{\tau=1}^{|\beta|} \left(\dfrac{-i\partial_{t}^{\ell_\tau} \mathscr{C}_{\xi}(t)}{\ell_\tau!}\right)^{k_\tau},
	\end{equation}
	where, for each integer  $1\leq k\leq |\beta|$,
	the sum is taken on the   set $S(\beta,k)$ consisting of elements
	\begin{equation*}
		(k_1,\ldots, k_{|\beta|}; \ell_1, \ldots , \ell_{|\beta|})\in \mathbb{Z}_+^{|\beta|}\times (\mathbb{Z}_+^n)^{|\beta|},
	\end{equation*}
	such that $\sum_{j=1}^{|\beta|}k_j=k$ and  $\sum_{\tau=1}^{|\beta|}k_\tau|\ell_\tau|=|\beta|$. Additional conditions on these elements are precisely described in Section 4 in  \cite{Silva2020}.
	
	\bigskip

	Note that,  $\partial_{t_j}\mathscr{C}_{\xi}(t) =(c_j(t)-c_{j0})p_j(\xi)$ and, in view of $c(t,\xi)$ being closed, we have  $p_j(\xi)\partial_{t_k}c_j(t)=p_k(\xi)\partial_{t_j}c_k(t)$, for all $j,k=1,\ldots,n$. Also, since each $\ell_\tau \in \mathbb{Z}_+^n$ in the previous formula is nonzero, $c_j$ are bounded functions and $p_j(\xi)$ satisfies \eqref{symbol pj},  then  there exists a constant $C_\tau>0$ such that
	$$ |\partial_{t}^{\ell_\tau} \mathscr{C}_{\xi}(t)| \leq C_{\tau}|\xi|^m,\quad \forall\xi \in \mathbb{Z}^N,
	$$
	where $m=\max\{0,m_1,\ldots,m_n\}$.
	
	Additionally, condition \eqref{condition D} implies 
	$$
	| e^{-i \mathscr{C}_{\xi}(t)}|=|e^{\Im \mathscr{C}_{\xi}(t)}| \leq C
	|\xi|^{\kappa}, \quad \forall\xi \in \mathbb{Z}^N.
	$$ 
	Thus, it follows from Formula  \eqref{Faa Di Bruno Formula} and the previous inequalities that there exists a constant $C_\beta>0$ such that
	\begin{equation}\label{inequality derivative exponential}
		|\partial_{t}^{\beta} e^{-i \mathscr{C}_{\xi}(t)}| \leq C_{\beta}
		|\xi|^{m |\beta| + \kappa}, \quad \forall\xi \in \mathbb{Z}^N.
	\end{equation}

	Given $M\in \mathbb{N}$, we take $N_0\in \mathbb{N}$ so that
	$m|\beta| + M < N_0$. Thus, since $u_K$ is a smooth function, there exists $C_0>0$ such that
	\begin{equation}\label{inequality derivative fourier coef}
		|\partial_{t}^{\alpha-\beta}\widehat{u_K}(t,\xi)| \leq C_0 |\xi|^{-N_0}, \quad \forall\xi \in \mathbb{Z}^{N}.
	\end{equation} 
	
	Therefore, using the inequalities \eqref{inequality derivative exponential} and \eqref{inequality derivative fourier coef} in formula  \eqref{Leibniz rule}, we may conclude the validity of
	\eqref{estimative derivatives}. Hence, $\Psi u\in C^{\infty}(\mathbb{T}^{n+N}, \wedge^{p,0})$.

	\bigskip
	
	On the other hand, consider
	$$
	u=\sum_{\xi\in \mathbb{Z}^N}
	\sum_{|K|=p}
	\left(
	\widehat{u_K}(t,\xi)e^{i\xi \cdot x}
	\right)dt_K \in \mathcal{D}'(\mathbb{T}^{n+N}, \wedge^{p,0}).
	$$

Since $u_K\in\mathcal{D}'(\mathbb{T}^{n+N})$, there exist $C>0$ and $M\in\mathbb{N}$ such that, for any $\varphi \in \mathcal{C}^{\infty}(\mathbb{T}_t^n)$
	we have
	$$
	|\langle e^{-i \mathscr{C}_{\xi}(t)}  \widehat{u_K}(t,\xi) \, , \,
	\varphi(t)\rangle | = |\langle \widehat{u_K}(t,\xi) \, , \,
	e^{-i \mathscr{C}_{\xi}(t)} \varphi(t)\rangle |  \leq C \rho_M(e^{-i \mathscr{C}_{\xi}} \varphi)|\xi|^{M},
	$$
	where
	$$
	\rho_M(e^{-i \mathscr{C}_{\xi}} \varphi) =  \sup_{|\alpha| \leq M}\left\{ \sup_{t\in \mathbb{T}^n} |\partial_t^\alpha(e^{-i \mathscr{C}_{\xi}(t)} \varphi(t))|\right\}.
	$$
	
	By similar arguments to those used in the derivatives  $\partial_t^\alpha(e^{-i \mathscr{C}_{\xi}(t)} \varphi(t))$, it is not difficult to see that
	$$
	\rho_M(e^{-i \mathscr{C}_{\xi}} \varphi) \leq C|\xi|^{mM+\kappa} \rho_M(\varphi).
	$$
	Then, for $M_0\geq mM+\kappa$ we get
	$$
	|\langle e^{-i \mathscr{C}_{\xi}(t)}  \widehat{u_K}(t,\xi) \, , \,
	\varphi(t)\rangle |  \leq C \rho_{M_0}(\varphi)|\xi|^{M_0},
	$$
	which implies that $\Psi u \in \mathcal{D}'(\mathbb{T}^{n+N}, \wedge^{p,0})$. Similar arguments can be applied to the inverse operator $\Psi^{-1}$, concluding that $\Psi$ is an isomorphism on the mentioned spaces.

	\bigskip
	
	Now, we verify the validity of  Formula \eqref{conjugation}. Consider 
	$u \in \mathcal{D}'(\mathbb{T}^{n+N}, \wedge^{p,0})$, or $u \in C^{\infty}(\mathbb{T}^{n+N}, \wedge^{p,0})$. Then, 
	\begin{align*}
		\mathbb{L}_0^p\circ \Psi^{-1}(u) & = \mathbb{L}_0^p (\Psi^{-1} u) =
		\sum_{\xi \in \mathbb{Z}^N}(d_t + ic_0(D_x)\wedge)(\widehat{u}(t,\xi)e^{i \mathscr{C}_{\xi}(t)} e^{i x \cdot \xi})\\
		& = \sum_{\xi \in \mathbb{Z}^N} e^{i \mathscr{C}_{\xi}(t)} 
		\left(d_t\widehat{u}(t,\xi) + i (c_{\xi0} + d_t\mathscr{C}_{\xi}(t))\right)e^{i x \cdot \xi}\\
		& = \Psi^{-1}  \mathbb{L}^p(u).
	\end{align*}

	Finally, we show that the spaces
	$\mathbb{E}^{p}_{0}$ and $\mathbb{E}^{p}$ are  isomorphic: if $f \in \mathbb{E}^{p}_{0}$, then  
	$$
	\mathbb{L}^{p+1} (\Psi f) = \mathbb{L}^{p+1} \circ \Psi (f) = 
	\Psi \circ \mathbb{L}^{p+1}_0(f)  = \Psi(\mathbb{L}^{p+1}_0 f) =0.
	$$
	Moreover, 
	$$
	\widehat{\Psi f}(t,\xi) e^{i(\psi_{\xi}(t) + \mathscr{C}_{\xi}(t))} = e^{-i\mathscr{C}_{\xi}(t)}
	f(t,\xi)e^{i(\psi_{\xi}(t) + \mathscr{C}_{\xi}(t))} = 
	f(t,\xi)e^{i\psi_{\xi}(t)}.
	$$
	
	Therefore, if $c_{\xi0}$ is integral, then $f(t,\xi)e^{i\psi_{\xi}(t)}$ is a  ($p+1$)-exact form, which implies the  same for  $\widehat{\Psi f}(t,\xi) e^{i(\psi_{\xi}(t) + \mathscr{C}_{\xi}(t))}$. Thus, we have $\Psi(\mathbb{E}^{p}_{0}) \subset \mathbb{E}^{p}$.

	With similar arguments we may obtain  $\Psi^{-1}(\mathbb{E}^{p}) \subset \mathbb{E}^{p}_{0}$, which concludes the proof.

\end{proof}

As an illustration of the result of Theorem \ref{reduction normal form thm}, we recall the reduction made in \cite{BerPet99jde}.

\begin{example}
	For $p=0,1,\ldots,n-1$, consider the differential operator	
	$$
	\mathbb{L}^p=d_t+a(t)\wedge \dfrac{\partial}{\partial x},
	$$	
	defined on $\mathbb{T}_t^{n}\times \mathbb{T}_x^{1}$, where $a(t)$ is a smooth, real and closed 1-form. 
	In this case, writing $a(t)=a_0+d_t\mathscr{A}(t)$, where $a_0=\sum_{j=1}^{n}a_{j0}dt_j$ is constant and $\mathscr{A}(t)$ is a real valued function, clearly we obtain  
	$\exp \left(\xi \Im \mathscr{A}(t) \right)   \equiv 1$, for all $\xi\in\mathbb{Z}$. Therefore,  $\mathbb{L}^p$ is globally solvable if and only if $$
	\mathbb{L}_0^p=d_t+a_0\wedge \dfrac{\partial}{\partial x}
	$$ is globally solvable. Thus, it follows from Theorem \ref{differential_complex} that $\mathbb{L}^p$ is globally solvable if and only if $a_0$ is either rational or non-Liouville.

\end{example}

%====================================
%====================================
\subsection{Decoupled systems \label{sec_dec_system}}
%====================================
%====================================

	In this section, we apply the results of Theorem \ref{reduction normal form thm} to study the global sol\-va\-bi\-lity of  
	$\mathbb{L}^p$ by assuming  that in System \eqref{main system}
	each coefficient   $c_j$ is a smooth function depending only on the variable $t_j \in \mathbb{T}_{t_j}^{1}$. Precisely, we consider the system
	\begin{equation*}\label{system uncloupled}
		L_j = D_{t_j} + c_j(t_j) P_j(D_x),\quad   j=1,\ldots,n.
	\end{equation*}

	In the sequel, we will use the following notations: 
	$$
	c_j(t_j) = a_j(t_j) + ib_j(t_j), \quad  p_j(\xi) = \alpha_j(\xi) + i \beta_j(\xi),\quad  \xi\in\mathbb{Z}^N,
	$$
	where $a_j, b_j$ are real valued functions and $\alpha_j, \beta_j$ are real valued sequences. Also, we set 
	$$
	\mathscr{A}_\xi(t) = \sum_{j=1}^{n} p_j(\xi) A_j(t_j) \ \textrm{ and } \
	\mathscr{B}_\xi(t) = \sum_{j=1}^{n} p_j(\xi) B_j(t_j),
	$$
	where
	$$
	A_j(t_j)  = \int_{0}^{t_j}a_j(s)ds - a_{j0} t_j \ \textrm{ and } \
	B_j(t_j)  = \int_{0}^{t_j}b_j(s)ds - b_{j0} t_j,
	$$
	and $c_{j0}=a_{j0}+ib_{j0}$, with $a_{j0},b_{j0}\in\mathbb{R}$.

	Therefore, for each $\xi\in\mathbb{Z}^N$, writing the closed 1-form  as
	$$
	c(t,\xi)=\sum_{j=1}^{n}c_j(t_j)p_j(\xi)dt_j=\sum_{j=1}^{n}\mathcal{M}_j(t_j,\xi)dt_j,
	$$
	we obtain $c(t,\xi) =c_{\xi0} + d_t\mathscr{C}_\xi(t)
	$, where 
	$$c_{\xi0} = \sum_{j=1}^{n}c_{j0}p_j(\xi)dt_j\text{\; and \;}\mathscr{C}_\xi(t) = \mathscr{A}_\xi(t) + i \mathscr{B}_\xi(t).
	$$ 
	
	We point out that condition $\mathscr{D}$, stated in Definition \ref{definition condition D}, is fulfilled if and only if for each $j\in\{1,\ldots,n\}$, there exist positive cons\-tants $C_j,\kappa_j$ and $n_0$ satisfying
	\begin{equation}\label{general-cond-reduc decoupled}
		\sup_{\zeta \in \mathbb{T}}\left\{ 
		\exp \left(\int_{0}^{\zeta} \Im \mathcal{M}_{j}(s,\xi) ds\right) \right\} \leq 
		C_j |\xi|^{\kappa_j}, 
	\end{equation}
	for all $|\xi|\geq n_0$.

	\bigskip

	Now,  let us display a class of systems  satisfying Condition  \eqref{general-cond-reduc decoupled}. For this, we recall  Definition 4.1 in \cite{AviGonKirMed19}: We say that a  function $\phi:\mathbb{Z}^N\rightarrow \mathbb{C}$ has at most logarithmic growth,
	if there are positive constants $C$ and $n_0$ such that
	\begin{equation}\label{log-estimate-1}
		|\phi(\xi)| \leq  C \log(|\xi|), \quad  |\xi|\geq  n_0.
	\end{equation}
	By simplicity, when this property is verified we write $\phi(\xi) = \mathcal{O}(\log(|\xi|))$. When it  is not, we say that $\phi(\xi)$ has super-logarithmic growth.

\begin{definition}\label{definition log}
	Consider the functions
	$$
	\Im\mathcal{M}_j(t_j,\xi)=a_j(t_j)\beta_j(\xi)+b_j(t_j)\alpha_j(\xi).
	$$
	
	We denote by $\mathcal{L}$ the subset of indexes $j\in\{1,\ldots,n\}$  such that  one of the following conditions holds:
	\begin{enumerate}
		\item [(i)] $p_j(\xi) = \mathcal{O}(\log(|\xi|))$;
		
		\item [(ii)] $\alpha_j(\xi) = \mathcal{O}(\log(|\xi|))$,  $\beta_j(\xi)$ has super-logarithmic growth and  $a_j(\cdot)$ does not change sign;
		
		\item [(iii)] $\alpha_j(\xi)$ has super-logarithmic growth, $\beta_j(\xi) = \mathcal{O}(\log(|\xi|))$ and $b_j(\cdot)$ does not change sign.
	\end{enumerate}
\end{definition}

The proof of Proposition 6 in \cite{AviMed21} shows that if there exists $j$ such that $p_j(\xi)=\mathcal{O}(\log(|\xi|))$,  
then Condition \eqref{general-cond-reduc decoupled} holds. Moreover, if there exists $j$ such that $(i)$ or $(ii)$ holds, then, by applying the same arguments to those  used in Proposition 7 in \cite{AviMed21}, we also obtain Condition \eqref{general-cond-reduc decoupled}. Therefore, we can state the following result:

\begin{proposition}\label{proposition reduction}
	Suppose that $\mathcal{L}= \{1,\ldots,n\}$. Then, $\mathbb{L}^p$ is globally solvable if and only if $\mathbb{L}_0^p$ (given in \eqref{normal form operator}) is globally solvable.
\end{proposition}

\begin{example}\label{exe-qui-slog-hor} 
	Let $p=0,1,2$. 	Consider   the operator $\mathbb{L}^p$  associated with the system
	$$
	L_j = D_{t_j} +  c_j(t_j) P_j(D_x), \quad j=1,2,3,
	$$
	defined on $\mathbb{T}^4$, where the coefficients $c_j$ and the toroidal symbols $p_j(\xi)$, $\xi\in\mathbb{Z}$,  are as follows:
	\begin{itemize}
		\item $p_1(\xi)= \log(1+ |\xi|)$;
		
		\item $p_2(\xi)= 1 +  i \xi$, $\Im c_3(t_3) \equiv 0$ and  $\Re c_3(t_3)>0$;
		
		\item $p_3(\xi) = \xi$, $\Im c_1(t_1) \equiv 0$.
		
	\end{itemize}

	Note that,  $p_1(\xi)=\mathcal{O}(\log(|\xi|))$ and $p_2$ and $p_3$ satisfy conditions $(ii)$ and $(iii)$ of Definition \ref{definition log}, respectively. Additionally, since $L_2$ is globally solvable in $\mathbb{T}^2$, it follows that the normal form $\mathbb{L}_0^p$ is globally solvable, as observed in Remark \ref{remark coef constants}. We conclude from Proposition \ref{proposition reduction} that $\mathbb{L}^p$ is globally solvable.

\end{example}

The next example shows that it is possible to achieve the reduction to normal form  considering other classes of pseudo-differential operators that do not satisfy the conditions presented in Definition \ref{definition log}.

\begin{example}
	Let $P(D_x)$ be a pseudo-differential operator  on $\mathbb{T}^1$ with symbol $p(\xi) =\alpha(\xi) + i\beta(\xi)$ defined as follows:
	$$
	\alpha(\xi) = 
	\begin{cases}
		\xi^{-1}, & \textrm{ if } \ \xi<0  \ \textrm{ is odd,} \\
		|\xi|, & \textrm{ if } \ \xi \leq 0  \ \textrm{ is even,} \\
		0, & \textrm{ if  } \ \xi > 0,
	\end{cases}
	\ \textrm{and } \
	\beta(\xi) = 
	\begin{cases}
		1, & \textrm{ if  } \ \xi \leq 0, \\
		\xi, & \textrm{ if } \ \xi > 0.
	\end{cases}
	$$

	Now, let $a(\cdot)$ and $b(\cdot)$ be two real-valued, smooth  and  negative functions on $\mathbb{T}^1$ with disjoint supports. Defining  $L= D_{s} + (a(s)+ ib(s))P(D_x)$, we obtain
$$
\Im \mathcal{M}(s,\xi) =
\begin{cases}
	a(s) + {b(s)}/{\xi}, & \textrm{ if } \ \xi<0   \textrm{ is odd,} \\
	a(s) + b(s)|\xi|, & \textrm{ if } \ \xi \leq 0  \ \textrm{ is even,} \\
	a(s) \xi, &  \textrm{ if } \ \xi > 0.
\end{cases}
$$

	Although none of the conditions $(i)$, $(ii)$, or $(iii)$ of Definition \ref{definition log} are fulfilled, we have $\Im \mathcal{M}(s,\xi)\leq 0$ and therefore \eqref{general-cond-reduc decoupled} holds.

\end{example}

\textbf{Acknowledgments}
The first author thanks the support provided by the National Council for Scientific and Technological Development - CNPq, Brazil (grants 423458/2021-3, 402159/2022-5, 200295/2023-3, and 305630/2022-9).

\bibliographystyle{plain} 
\bibliography{references}

\end{document}